\documentclass[a4paper,11pt]{article}
\usepackage{microtype}
\usepackage[english,italian]{babel}
\usepackage[dvips]{graphicx}
\usepackage{float}
\usepackage{amsfonts}
\usepackage{amssymb}
\usepackage{amsmath}
\usepackage{amsthm}
\usepackage{enumerate}
\usepackage{MnSymbol}

\newcommand{\G}{\Gamma}

\newcommand{\eps}{\varepsilon}
\newcommand{\cOx}{{\cal O}_x}
\newcommand{\cOy}{{\cal O}_y}
\newcommand{\cEx}{{\cal E}_x}
\newcommand{\cEy}{{\cal E}_y}

\title{On $d$-divisible graceful $\alpha$-labelings of $C_{4k}\times P_m$}

\author{Anita Pasotti \thanks{Dipartimento di Matematica,
Facolt\`a di Ingegneria, Universit\`a degli Studi di Brescia,
Via Valotti, 9, I-25133 Brescia, Italy. email: anita.pasotti@ing.unibs.it}}

\date{}
\newtheorem{defi}{Definition}[section]
\newtheorem{prop}[defi]{Proposition}

\newtheorem{ex}[defi]{Example}
\newtheorem{thm}[defi]{Theorem}

\begin{document}
\selectlanguage{english}
\maketitle
\selectlanguage{english}

\begin{abstract}
In \cite{APArs} the concept of a $d$-divisible graceful $\alpha$-labeling has been introduced
as a generalization of classical $\alpha$-labelings
and it has been shown how it is useful to obtain certain cyclic graph decompositions.
In the present paper it is proved the existence of $d$-divisible graceful $\alpha$-labelings
of $C_{4k}\times P_m$ for any integers $k\geq1$, $m\geq2$ for several values of $d$.
\end{abstract}


\section{Introduction}
We assume familiarity with the basic concepts about graphs.\\
As usual, we denote by $K_v$ and $K_{m\times n}$ the \emph{complete
graph on $v$ vertices} and the \emph{complete $m$-partite graph with parts of size $n$},
respectively.
Also, let $C_k$, $k\geq 3$, be the cycle on $k$ vertices
and let $P_m$, $m \geq2$, be the path on $m$ vertices.
Graphs of the form $C_k \times P_m$ can be viewed as grids on cylinders
and they are bipartite if and only if $k$ is even. If $m=2$,
$C_k \times P_2$ is nothing but the prism $T_{2k}$ on $2k$ vertices.
For any graph $\G$ we write $V(\G)$ for the set of its vertices and $E(\G)$
for the set of its edges. If $|E(\G)|=e$, we say that $\G$ has \emph{size} $e$.

Given a subgraph $\G$ of a graph $K$, a $\G$-\emph{decomposition of} $K$
is a set of graphs, called \emph{blocks}, isomorphic to $\G$ whose edges
partition the edge-set of $K$. Such a decomposition is said to be \emph{cyclic} when it is invariant
under a cyclic permutation of all vertices of $K$. In the case that $K=K_v$
one also speaks of a $\G$-\emph{system of order} $v$. The problem of establishing
the set of values of $v$ for which such a system exists is in general quite difficult.
For a survey on graph decompositions see \cite{BE}.

The concept of a \emph{graceful labeling} of $\G$, introduced by A. Rosa \cite{R},
is quite related to the existence problem of cyclic $\G$-systems. A \emph{graceful
labeling} of a graph $\G$ of size $e$ is an injective function $f: V(\G)\rightarrow
\{0,1,2,\ldots,e\}$ such that
$$\{|f(x)-f(y)| \ |\ [x,y]\in E(\G)\}=\{1,2,\ldots,e\}.$$
In the case that $\G$ is bipartite and $f$ has the additional property that its
maximum value on one of the two bipartite sets does not reach its minimum on the other
one, one says that $f$ is an $\alpha$-labeling. In \cite{R}, Rosa proved that
if a graph $\G$ of size $e$ admits a graceful labeling then there exists a cyclic
$\G$-system of order $2e+1$ and that if it
admits an $\alpha$-labeling then there exists
a cyclic $\G$-system of order $2en+1$ for any positive integer $n$.
For a very rich survey on graceful labelings we refer to \cite{G}.

Many variations of graceful labelings have been considered. In particular Gnana
Jothi \cite{GJ} defines an \emph{odd graceful labeling} of a graph $\G$ of size $e$
as an injective function $f: V(\G)\rightarrow
\{0,1,2,\ldots,2e-1\}$ such that
$$\{|f(x)-f(y)| \ |\ [x,y]\in E(\G)\}=\{1,3,5,\ldots,2e-1\}.$$
In a recent paper, see \cite{APArs}, we have introduced the following new definition
which is, at the same time, a generalization of the concepts of a graceful labeling
(when $d=1$) and of an odd graceful labeling (when $d=e$).
\begin{defi}
Let $\G$ be a graph of size
$e=d\cdot m$. A $d$-\emph{divisible graceful labeling} of $\G$
is an injective function $f:V(\G) \rightarrow \{0,1,2,\ldots, d(m+1)-1\}$
such that
\begin{align*}
\{|f(x)-f(y)|\ |\ [x,y]\in E(\G)\} &=\{1,2,3,\ldots,d(m+1)\}\\
&\quad  - \{m+1,2(m+1),\ldots,d(m+1)\}.
\end{align*}
Namely the set $\{|f(x)-f(y)|\ |\ [x,y]\in E(\G)\} $
can be divided into $d$ parts $P^0,P^1,\ldots,P^{d-1}$
where $P^i:=\{(m+1)i+1,(m+1)i+2,\ldots,(m+1)i+m\}$ for any
$i=0,1,\ldots, d-1$.
\end{defi}
\noindent
The $\alpha$-labelings can be generalized in a similar way.
\begin{defi}
A $d$-\emph{divisible graceful} $\alpha$-\emph{labeling} of a bipartite graph $\G$ is a
$d$-divisible graceful labeling of $\G$  having the property that
its maximum value on one of the two bipartite sets does not reach its minimum value on the other one.
\end{defi}

\noindent
We have to point out that in \cite{APArs} the above labelings have been called
``$d$-graceful ($\alpha$-)labelings'', but the author was unaware that
this name is already used
in the literature whit a different meaning,
see \cite{MT} and \cite{S}.\\
\\
It is known that there is a close relationship
between graceful labelings and difference families, see \cite{AB}.
In \cite{APArs} we established relations between $d$-divisible graceful ($\alpha$-)labelings
and a generalization of difference families introduced in \cite{BP},
proving the following theorems.
\begin{thm}
If there exists a $d$-divisible graceful labeling of a graph $\G$ of size $e$ then there exists
a cyclic $\G$-decomposition of $K_{\left(\frac{e}{d}+1\right)\times 2d}$.
\end{thm}

\begin{thm}\label{decomp}
If there exists a $d$-divisible graceful $\alpha$-labeling of a graph $\G$ of size $e$ then there exists
a cyclic $\G$-decomposition of $K_{(\frac{e}{d}+1)\times 2dn}$ for any integer $n\geq1$.
\end{thm}
\noindent
In this paper we determine the existence of $d$-divisible graceful $\alpha$-labelings of
$C_{4k}\times P_m$ for several values of $d$. In order to obtain these results,
first of all we will find $d$-divisible graceful $\alpha$-labelings of prisms, which correspond to
the case $m=2$, and then by induction on $m$
we will be able to construct $d$-divisible graceful $\alpha$-labelings of $C_{4k}\times P_m$ for any
$m\geq 2$. For what said above, these results allow us to obtain new infinite classes of cyclic
decompositions of the complete multipartite graph in copies of $C_{4k}\times P_m$.

\section{$d$-divisible graceful $\alpha$-labelings of prisms}
In this section we will investigate the existence of $d$-divisible graceful $\alpha$-labelings
of prisms. From now on, given two integers $a$ and $b$, by $[a,b]$ we will denote the
set of integers $x$ such that $a\leq x \leq b$.\\
For convenience, we denote the $2k$ vertices of $T_{2k}$ by $x_1,x_2,\ldots,x_k$;
$y_1,y_2,\ldots,$ $y_k$ where the $x_i$'s are the consecutive vertices of one $k$-cycle
and the $y_i$'s are consecutive vertices of the other $k$-cycle and $x_i$ is
connected to $y_i$. Clearly $T_{2k}$ has
size $e=3k$ and it is bipartite if and only if $k$ is even.
In \cite{FG} Frucht and Gallian proved that $T_{2k}$ admits an $\alpha$-labeling
if and only if $k$ is even.

\begin{thm}\label{prism3}
The prism $T_{8k}$ admits a $3$-divisible graceful $\alpha$-labeling for every $k\geq1$.
\end{thm}
\noindent
Proof. We set
$\cOx=\{x_1,x_3,\ldots,x_{4k-1}\}$, $\cEx=\{x_2,x_4,\ldots,x_{4k}\}$,
$\cOy=\{y_1,y_3,$ $\ldots,y_{4k-1}\}$, $\cEy=\{y_2,y_4,\ldots,$ $y_{4k}\}$.
Clearly $\cOx\cup \cEy$ and $\cOy\cup \cEx$ are the two bipartite sets of $V(T_{8k})$.\\
Consider the map $f:V(T_{8k})\rightarrow\{0,1,\ldots, 12k+2\}$ defined as follows:
\begin{align*}
f(x_{2i+1})&=\left\{
\begin{array}{l}
6k+1 \\
8k+2-i\\
8k+1-i
\end{array}
\right. & &
\begin{array}{l}
\textrm{for}\ i=0\\
 \textrm{for}\ i \in \left[ 1,k \right]\\
  \textrm{for}\ i \in \left[ k+1,2k-1 \right]
\end{array}\\
f(x_{2i})&=4k+i & &
\begin{array}{l}
\textrm{for}\ i\in [1,2k].
\end{array}\\
f(y_{2i+1})&=i & &
\begin{array}{l}
\textrm{for}\ i\in [0,2k-1]
\end{array}\\
f(y_{2i})&=\left\{
\begin{array}{l}
12k+3-i\\
12k+2-i
\end{array}
\right. & &
\begin{array}{l}
 \textrm{for}\ i \in \left[ 1,k \right]\\
  \textrm{for}\ i \in \left[ k+1,2k \right].
\end{array}
\end{align*}
We have
\begin{align*}
f(\cOy\cup\cEx)&=[0,2k-1]\cup[4k+1,6k]\\
f(\cOx\cup\cEy)&=[6k+1,7k]\cup[7k+2,8k+1]\cup[10k+2,11k+1]\cup\\
&\cup[11k+3,12k+2].
\end{align*}
Hence $f$ is injective and $\max f(\cOy\cup\cEx)<\min f(\cOx\cup\cEy)$. Now for $i=1,\ldots,4k$ set
\begin{eqnarray}\label{ers}
\sigma_i=|f(x_{i+1})-f(x_i)|,\quad \eps_i=|f(y_{i+1})-f(y_i)|,\quad \rho_i=|f(x_{i})-f(y_i)|
\end{eqnarray}
where the indices are understood modulo $4k$.
By a direct calculation, one can see that
\begin{align*}
\sigma_1 &=2k,\\
\{\sigma_i\ |\ i=2,\ldots,2k+1\} &=[2k+1,4k]\\
\{\sigma_i\ |\ i=2k+2,\ldots,4k\} &=[1,2k-1]\\
\rho_1 &= 6k+1\\
\{\rho_i\ |\ i=2,\ldots,2k+1\} &=[6k+2,8k+1],\\
\{\rho_i\ |\ i=2k+2,\ldots,4k\} &=[4k+2,6k],\\
\{\eps_i\ |\ i=1,\ldots,2k\} &=[10k+3,12k+2],\\
\{\eps_i\ |\ i=2k+1,\ldots,4k-1\} & =[8k+3,10k+1],\\
\eps_{4k} & =10k+2.
\end{align*}
Hence
$\{\sigma_i\ |\ i=1,\ldots,4k\}=[1,4k]$,
 $\{\rho_i\ |\ i=1,\ldots,4k\}=[4k+2,8k+1]$ and
 $\{\eps_i\ |\ i=1,\ldots,4k\}=[8k+3,12k+2]$.
 This concludes the proof.
\hfill$\Box$

\begin{thm}\label{prism6}
The prism $T_{8k}$ admits a $6$-divisible graceful $\alpha$-labeling for every $k\geq1$.
\end{thm}
\noindent
Proof. Set $\cOx,\cEx,\cOy,\cEy$ as in the proof of previous theorem.
Consider the map
$f:V(T_{8k})\rightarrow \{0,1,\ldots,12k+5\}$ defined as follows:
\begin{align*}
f(x_{2i+1})&=\left\{
\begin{array}{l}
6k+2\\
8k+4-i\\
8k+2-i
\end{array}
\right. & &
\begin{array}{l}
\textrm{for}\ i=0\\
 \textrm{for}\ i \in \left[ 1,k \right]\\
  \textrm{for}\ i \in \left[ k+1,2k-1 \right]
\end{array}\\
f(x_{2i})&=4k+1+i & &
\begin{array}{l}
\textrm{for}\ i\in [1,2k].
\end{array}\\
f(y_{2i+1})&=i & &
\begin{array}{l}
\textrm{for}\ i\in [0,2k-1]
\end{array}\\
f(y_{2i})&=\left\{
\begin{array}{l}
12k+6-i\\
12k+4-i
\end{array}
\right. & &
\begin{array}{l}
 \textrm{for}\ i \in \left[ 1,k \right]\\
  \textrm{for}\ i \in \left[ k+1,2k \right].
\end{array}\\
\end{align*}
It results
\begin{align*}
f(\cOy\cup\cEx)&=[0,2k-1]\cup[4k+2,6k+1]\\
f(\cOx\cup\cEy)&=[6k+2,7k+1]\cup[7k+4,8k+3]\cup[10k+4,11k+3]\cup\\
&\cup[11k+6,12k+5].
\end{align*}
Hence $f$ is injective and $\max f(\cOy\cup\cEx)<\min f(\cOx\cup\cEy)$.
Let $\eps_i, \rho_i, \sigma_i$, for $i=1,\ldots,4k$, be as in (\ref{ers}).
It is not hard to see that
\begin{align*}
\{\sigma_i\ |\ i=1,\ldots,4k\} &=[1,2k]\cup[2k+2,4k+1]\\
\{\rho_i\ |\ i=1,\ldots,4k\} &=[4k+3,6k+2]\cup[6k+4,8k+3]\\
\{\eps_i\ |\ i=1,\ldots,4k\} &=[8k+5,10k+4]\cup[10k+6,12k+5].
\end{align*}
Hence $f$ is a $6$-divisible graceful $\alpha$-labeling of $T_{8k}$.
\hfill$\Box$

\begin{thm}\label{prism12}
The prism $T_{8k}$ admits a $12$-divisible graceful $\alpha$-labeling for every $k\geq1$.
\end{thm}
\noindent
Proof. Also here we set $\cOx,\cEx,\cOy,\cEy$ as in the proof of Theorem \ref{prism3}.
We are able to prove the existence of a $12$-divisible graceful $\alpha$-labeling of $T_{8k}$
by means of two direct constructions where we distinguish the two cases: $k$ even
and $k$ odd.\\
\\
Case 1: $k$ even.\\
Consider the map $f:V(T_{8k})\rightarrow \{0,1,\ldots,12k+11\}$ defined as follows:
\begin{align*}
f(x_{2i+1})&=\left\{
\begin{array}{l}
6k+5\\[2pt]
8k+8-i\\[2pt]
8k+7-i\\[2pt]
8k+5-i
\end{array}
\right. & &
\begin{array}{l}
\textrm{for}\ i=0\\[2pt]
 \textrm{for}\ i \in \left[ 1,\frac{k}{2} \right]\\[2pt]
\textrm{for}\ i \in \left[ \frac{k}{2}+1,k \right]\\[2pt]
  \textrm{for}\ i \in \left[k+1,2k-1 \right]
\end{array}\\
f(x_{2i})&=\left\{
\begin{array}{l}
4k+3+i\\[2pt]
4k+4+i
\end{array}
\right. & &
\begin{array}{l}
\textrm{for}\ i\in \left[1,\frac{3k}{2}\right]\\[2pt]
\textrm{for}\ i\in \left[\frac{3k}{2}+1,2k\right].
\end{array}\\
f(y_{2i+1})&=\left\{
\begin{array}{l}
i\\[2pt]
i+1
\end{array}
\right. & &
\begin{array}{l}
\textrm{for}\ i\in \left[0,\frac{3k}{2}-1\right]\\[2pt]
\textrm{for}\ i\in \left[\frac{3k}{2},2k-1\right]
\end{array}\\
f(y_{2i})&=\left\{
\begin{array}{l}
12k+12-i\\[2pt]
12k+11-i\\[2pt]
12k+9-i
\end{array}
\right. & &
\begin{array}{l}
 \textrm{for}\ i \in \left[ 1,\frac{k}{2} \right]\\[2pt]
 \textrm{for}\ i \in \left[ \frac{k}{2}+1,k \right]\\[2pt]
 \textrm{for}\ i \in \left[k+1,2k \right]
\end{array}
\end{align*}
It is easy to see that
\begin{align*}
f(\cOy) &=\left[0,\frac{3k}{2}-1\right] \cup \left[\frac{3k}{2}+1,2k\right]\\
f(\cEx) &=\left[4k+4,\frac{11k}{2}+3\right] \cup \left[\frac{11k}{2}+5,6k+4\right]\\
f(\cOx) &=\left[6k+5,7k+4\right] \cup \left[7k+7,\frac{15k}{2}+6\right] \cup \left[\frac{15k}{2}+8,8k+7\right]\\
f(\cEy) &=\left[10k+9,11k+8\right] \cup \left[11k+11,\frac{23k}{2}+10\right] \cup \left[\frac{23k}{2}+12,12k+11\right].
\end{align*}
Hence $f$ is injective and $\max f(\cOy\cup\cEx)=6k+4<6k+5=\min f(\cOx\cup\cEy)$.
Set $\sigma_i, \eps_i, \rho_i$, for $i=1,\ldots,4k$, as in (\ref{ers}).
By a long and tedious calculation, one can see that
\begin{align*}
\{\sigma_i\ |\ i=1,\ldots,4k\} &=[1,4k+3]-\{k+1,2k+2,3k+3\}\\
\{\rho_i\ |\ i=1,\ldots,4k\} &=[4k+5,8k+7]-\{5k+5,6k+6,7k+7\}\\
\{\eps_i\ |\ i=1,\ldots,4k\} &=[8k+9,12k+11]-\{9k+9,10k+10,11k+11\}.
\end{align*}
This concludes the proof of Case 1.\\
\\
Case 2: $k$ odd.\\
Let now $f:V(T_{8k})\rightarrow \{0,1,\ldots,12k+11\}$ defined as follows:
\begin{align*}
f(x_{2i+1})&=\left\{
\begin{array}{l}
6k+5\\[2pt]
8k+8-i\\[2pt]
8k+6-i\\[2pt]
8k+5-i
\end{array}
\right. & &
\begin{array}{l}
\textrm{for}\ i=0\\[2pt]
 \textrm{for}\ i \in \left[ 1,k \right]\\[2pt]
\textrm{for}\ i \in \left[ k+1,\frac{3k-1}{2} \right]\\[2pt]
  \textrm{for}\ i \in \left[\frac{3k+1}{2},2k-1 \right]
\end{array}\\
f(x_{2i})&=\left\{
\begin{array}{l}
4k+3+i\\[2pt]
4k+4+i
\end{array}
\right. & &
\begin{array}{l}
\textrm{for}\ i\in \left[1,\frac{k+1}{2}\right]\\[2pt]
\textrm{for}\ i\in \left[\frac{k+3}{2},2k\right].
\end{array}\\
f(y_{2i+1})&=\left\{
\begin{array}{l}
i\\[2pt]
i+1
\end{array}
\right. & &
\begin{array}{l}
\textrm{for}\ i\in \left[0,\frac{k-1}{2}\right]\\[2pt]
\textrm{for}\ i\in \left[\frac{k+1}{2},2k-1\right]
\end{array}\\
f(y_{2i})&=\left\{
\begin{array}{l}
12k+12-i\\[2pt]
12k+10-i\\[2pt]
12k+9-i
\end{array}
\right. & &
\begin{array}{l}
 \textrm{for}\ i \in \left[ 1,k \right]\\[2pt]
 \textrm{for}\ i \in \left[k+1,\frac{3k-1}{2} \right]\\[2pt]
 \textrm{for}\ i \in \left[\frac{3k+1}{2},2k \right]
\end{array}
\end{align*}
Arguing exactly as in Case 1, one can check that $f$ is a $12$-divisible graceful $\alpha$-labeling
of $T_{8k}$.
\hfill$\Box$\\

\begin{figure}[H]
\begin{center}
\includegraphics[width=0.27\textwidth]{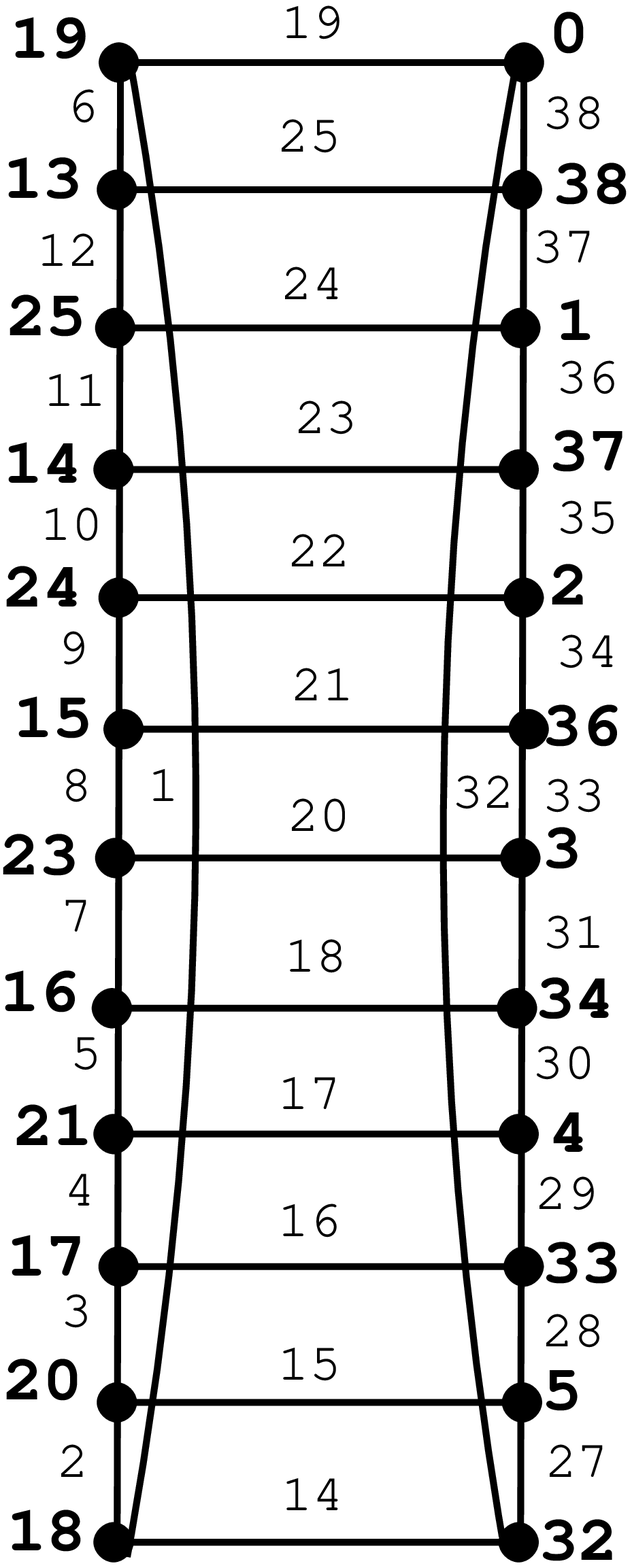}\quad\quad
\includegraphics[width=0.27\textwidth]{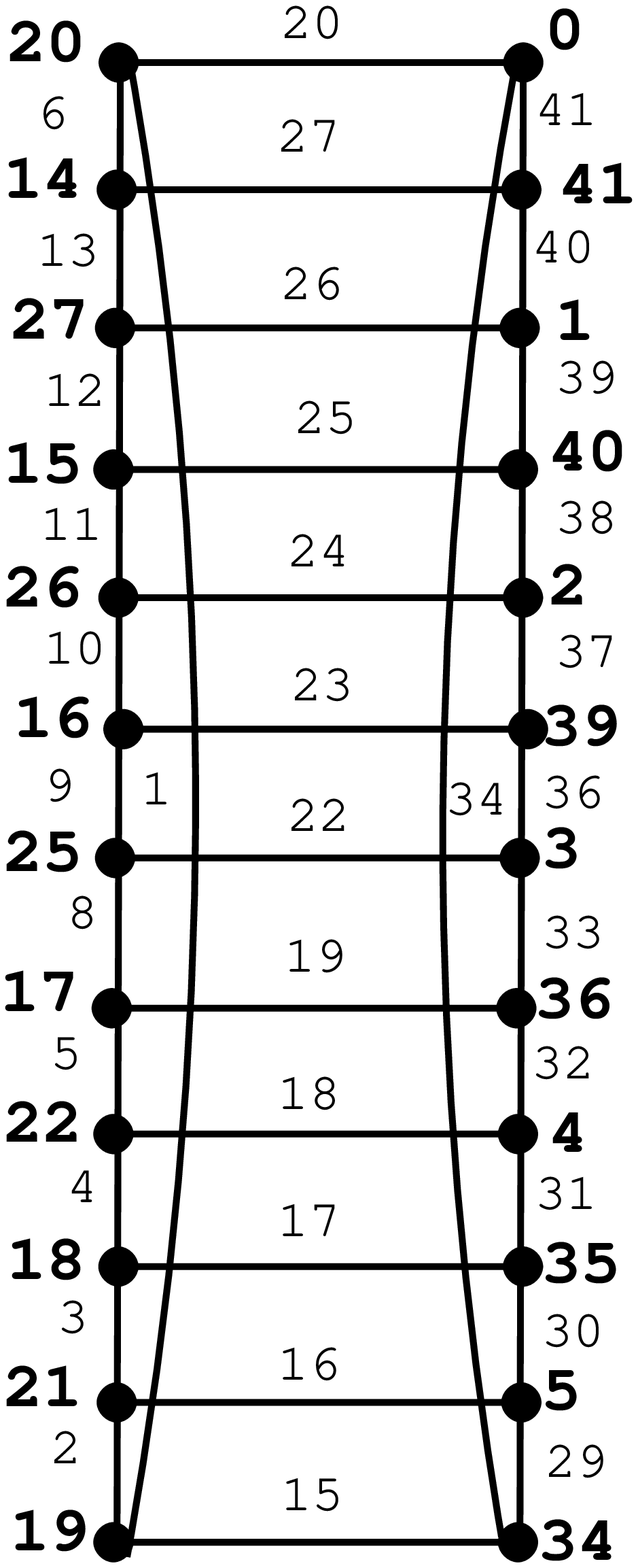}\quad\quad
\includegraphics[width=0.27\textwidth]{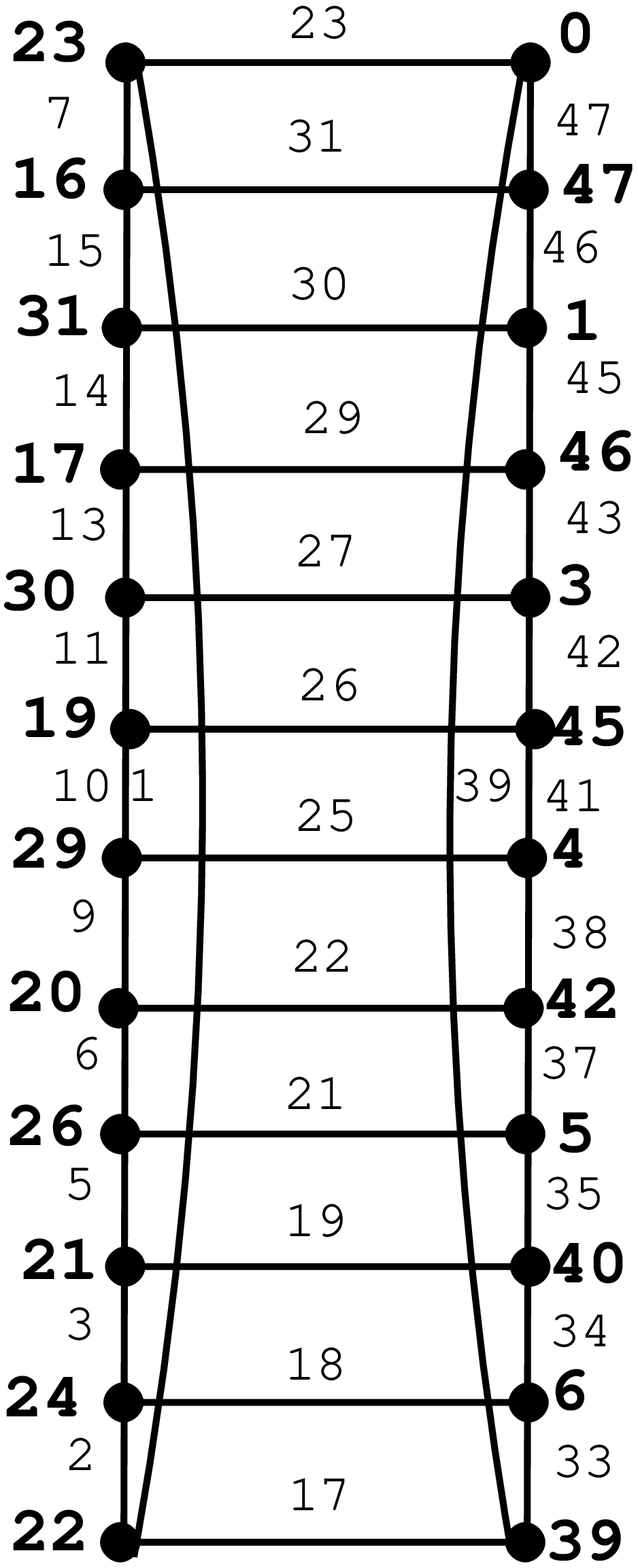}
\caption{$T_{24}$ }
\label{T24}
\end{center}
\end{figure}

\begin{ex}
The three graphs in Figure \ref{T24} show the $3$-divisible graceful $\alpha$-labeling,
the $6$-divisible graceful $\alpha$-labeling and the $12$-divisible graceful $\alpha$-labeling
of $T_{24}$ provided by previous theorems.
\end{ex}

\section{$d$-divisible graceful $\alpha$-labelings of $C_{4k}\times P_m$}
\begin{figure}[H]
\begin{center}
\includegraphics[width=0.55\textwidth]{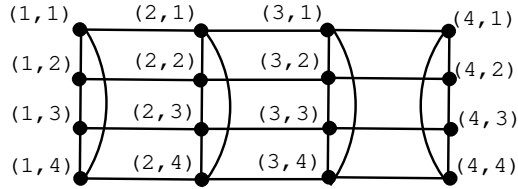}
\caption{$C_{4}\times P_4$ }
\label{C4P4ij}
\end{center}
\end{figure}

In this section using the results of the previous one we will construct
$d$-divisible graceful $\alpha$-labelings of $C_{4k}\times P_m$. In particular, since $e=4k(2m-1)$
we consider $d=2m-1,\ 2(2m-1),\ 4(2m-1)$. In \cite{JR} Jungreis and Reid
proved that for any $k,m\geq2$ not both odd there exists an $\alpha$-labeling
of $C_{2k}\times P_m$.\\
For convenience, we denote the vertices of $C_{4k}\times P_m$ as illustrated in Figure
\ref{C4P4ij} and we set
$C^i=((i,1),(i,2),\ldots,(i,4k))$ for any $i=1,\ldots,m$.

\begin{thm}\label{d=2m-1}
For any integer $k\geq1$ and $m\geq 2$, $C_{4k}\times P_m$ admits a
$(2m-1)$-divisible graceful $\alpha$-labeling.
\end{thm}
\noindent
Proof. We will prove the result by induction on $m$.
If $m=2$ the thesis follows from Theorem \ref{prism3}.
Let now $m\geq2$. Suppose that there exists a $(2m-1)$-divisible graceful $\alpha$-labeling
$f$ of $C_{4k}\times P_m$ with vertices of $C^m$ so labeled:
\begin{align*}
f(C^m) = &(0,(4k+1)(2m-1)-1,1,(4k+1)(2m-1)-2,2,\ldots,\\
&(4k+1)(2m-1)-k,k,(4k+1)(2m-1)-(k+2),
k+1,\dots,\\ &2k-1,(4k+1)(2m-1)-(2k+1)).
\end{align*}
Note that the $3$-divisible graceful $\alpha$-labeling of $C_{4k}\times P_2$ constructed in
Theorem \ref{prism3} has this property, in fact
$f(C^2)=(0,12k+2,1,12k+1,2,\ldots,k-1,11k+3,k,11k+1,k+1,\ldots,2k-1,10k+2)$.
So in order to obtain the thesis it is sufficient
to construct a $(2m+1)$-divisible graceful
$\alpha$-labeling $g$ of $C_{4k}\times P_{m+1}$ satisfying the
same property, namely such that
\begin{eqnarray}
\nonumber g(C^{m+1})&= &(0,(4k+1)(2m+1)-1,1,(4k+1)(2m+1)-2,2,\ldots,\\
\nonumber &\ &(4k+1)(2m+1)-k,k,
(4k+1)(2m+1)-(k+2),k+1,\ldots,\\
 &\ &2k-1,(4k+1)(2m+1)-(2k+1)).\label{gC(m+1)}
\end{eqnarray}
We set
$$g((i,j))=f((i,j))+(4k+1)\quad \forall i=1,\ldots,m,\ \forall j=1,\ldots,4k.$$
By the hypothesis on $f(C^m)$ it results
\begin{align*}
g(C^m)=&(4k+1,(4k+1)2m-1,4k+2,(4k+1)2m-2,4k+3,\ldots,\\
& (4k+1)2m-k, 5k+1,(4k+1)2m-(k+2),5k+2,\ldots,\\
& 6k,(4k+1)2m-(2k+1)).
\end{align*}
So there exists $j\in[1,4n]$ such that $g((m,j))=(4k+1)2m-1$.
We set $g(C^{m+1})$ as in (\ref{gC(m+1)}) where $g((m+1,j))=0$.\\
Now we will see that $g:V(C_{4k}\times P_{m+1})\rightarrow \{0,\ldots,(4k+1)(2m+1)-1\}$
defined as above
is indeed a $(2m+1)$-divisible graceful $\alpha$-labeling of $C_{4k}\times P_{m+1}$.
Since $f(V(C_{4k}\times P_m))\subseteq[0,(4k+1)(2m-1)-1]$, by the definition of $g$,
it follows that
$$g(V(C^1\cup C^2\cup\ldots\cup C^m))\subseteq [4k+1,(4k+1)2m-1].$$
Also we have
$$g(V(C^{m+1}))\subseteq[0,2k-1]\cup[(4k+1)(2m+1)-(2k+1),(4k+1)(2m+1)-1].$$
hence $g$ is an injective function.
Since, by hypothesis $f$ is a $(2m-1)$-divisible graceful $\alpha$-labeling of $C_{4k}\times P_m$, we have
$V(C_{4k}\times P_{m})=A \cupdot B$ with $\max_A f< \min_B f$.
Let $V(C_{4k}\times P_{m+1})=C \cupdot D$. By the construction, it follows that
\begin{align*}
g(C)&=(f(A)+(4k+1))\cup[0,2k-1]\\
g(D)&\subseteq(f(B)+(4k+1))\cup[(4k+1)(2m+1)-(2k+1),(4k+1)(2m+1)-1]
\end{align*}
hence $\max_C g< \min_D g$. Now we have to consider the differences between adjacent vertices.
Since $f$ is a $(2m-1)$-divisible graceful $\alpha$-labeling of $C_{4k}\times P_m$, by the construction of $g$, it results
$$\bigcup_{\begin{subarray}{c} i\in[1,m-1]\\ j\in[1,4k]\end{subarray}}|f((i,j))-f((i+1,j))|\ \cup
\bigcup_{\begin{subarray}{c} i\in[1,m]\\ j\in[1,4k]\end{subarray}} |f((i,j))-f((i,j+1))|=$$
$$[1,(4k+1)(2m-1)]-\{\beta(4k+1)\ |\ \beta\in[1,2m-1]\}$$
where the index $j$ is taken modulo $4k$.
Finally it is not hard to check that
\begin{align*}
&\{|g((m,j))-g((m+1,j))|\ |\ j\in[1,4k]\}=\\
&[(4k+1)2m-4k,(4k+1)2m-1]
\end{align*}
and
\begin{align*}
&\{|g((m+1,j))-g((m+1,j+1))|\ |\ j\in[1,4k]\}=\\
&[(4k+1)(2m+1)-4k,(4k+1)(2m+1)-1].
\end{align*}
This concludes the proof.
\hfill$\Box$

\begin{ex}
In Figure \ref{C4} we will show the $5$-divisible graceful $\alpha$-labeling of $C_4\times P_3$, the
$7$-divisible graceful $\alpha$-labeling of $C_4\times P_4$ and the
$9$-divisible graceful $\alpha$-labeling of $C_4\times P_5$ obtained
starting from the $3$-divisible graceful $\alpha$-labeling of $T_8=C_4\times P_2$
and following the construction illustrated in the proof
of Theorem $\ref{d=2m-1}$.
\begin{figure}[H]
\begin{center}
\includegraphics[width=0.3\textwidth]{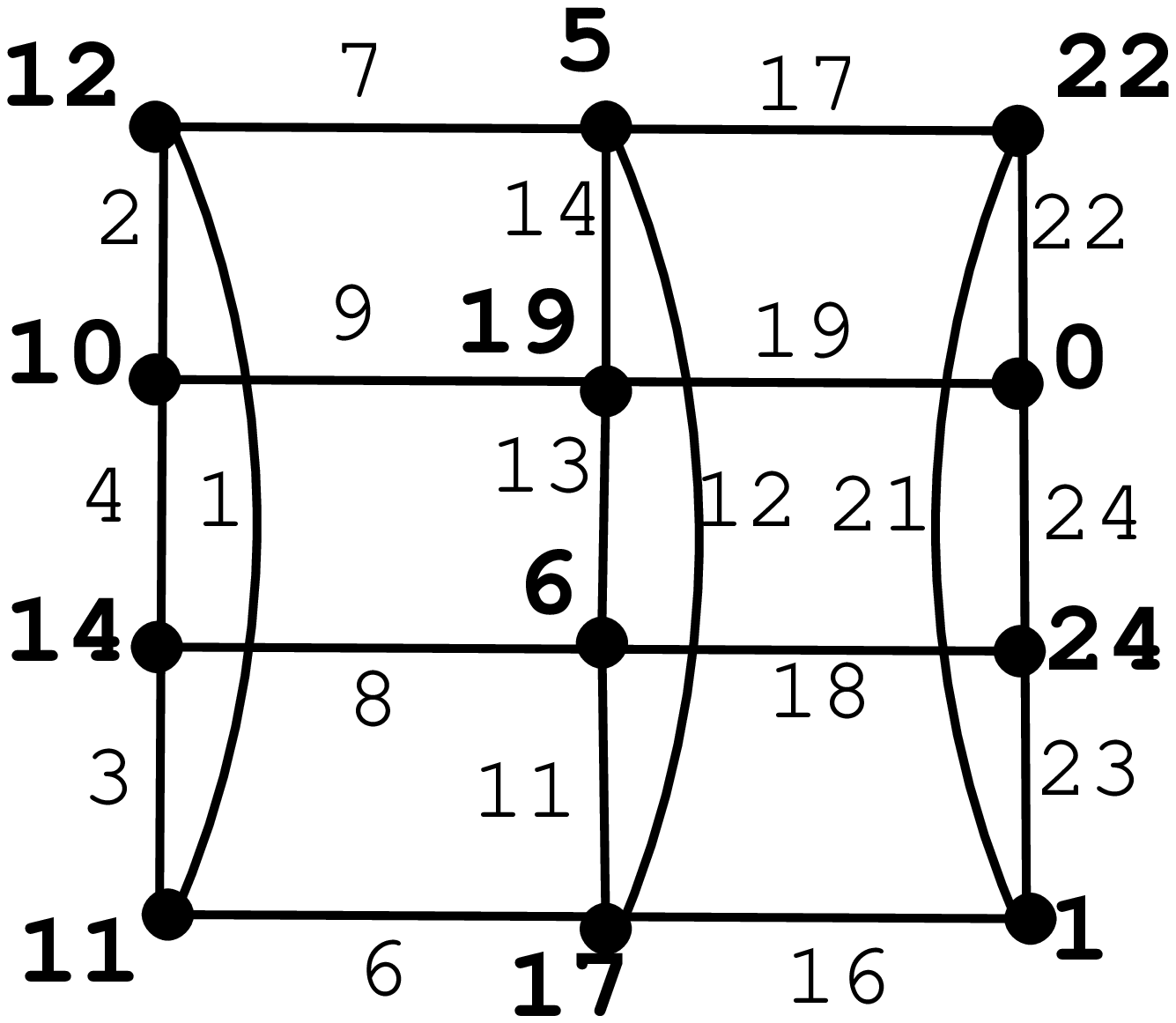}\quad\quad\quad
\includegraphics[width=0.4\textwidth]{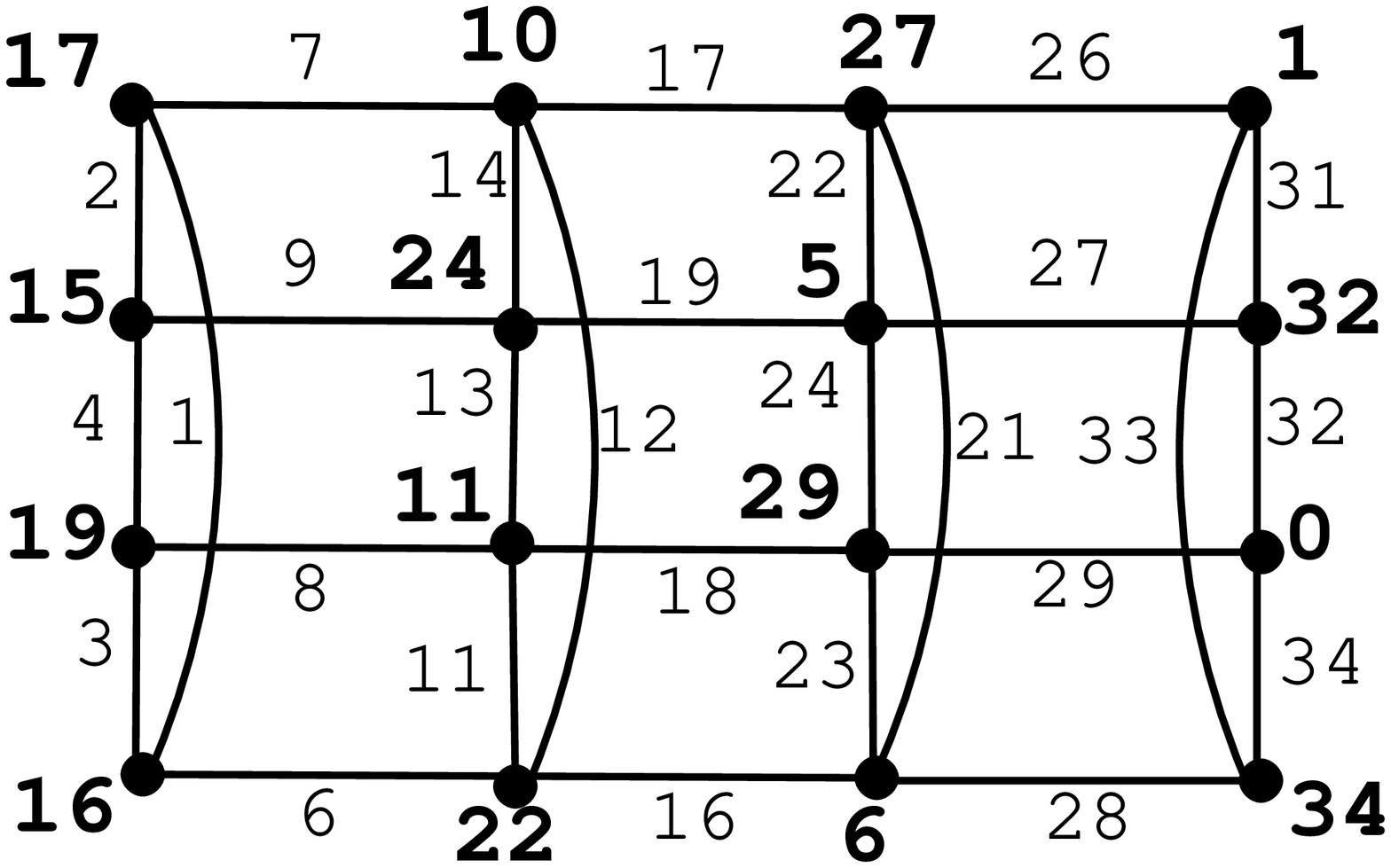}\\
\vspace{0.3cm}
\includegraphics[width=0.5\textwidth]{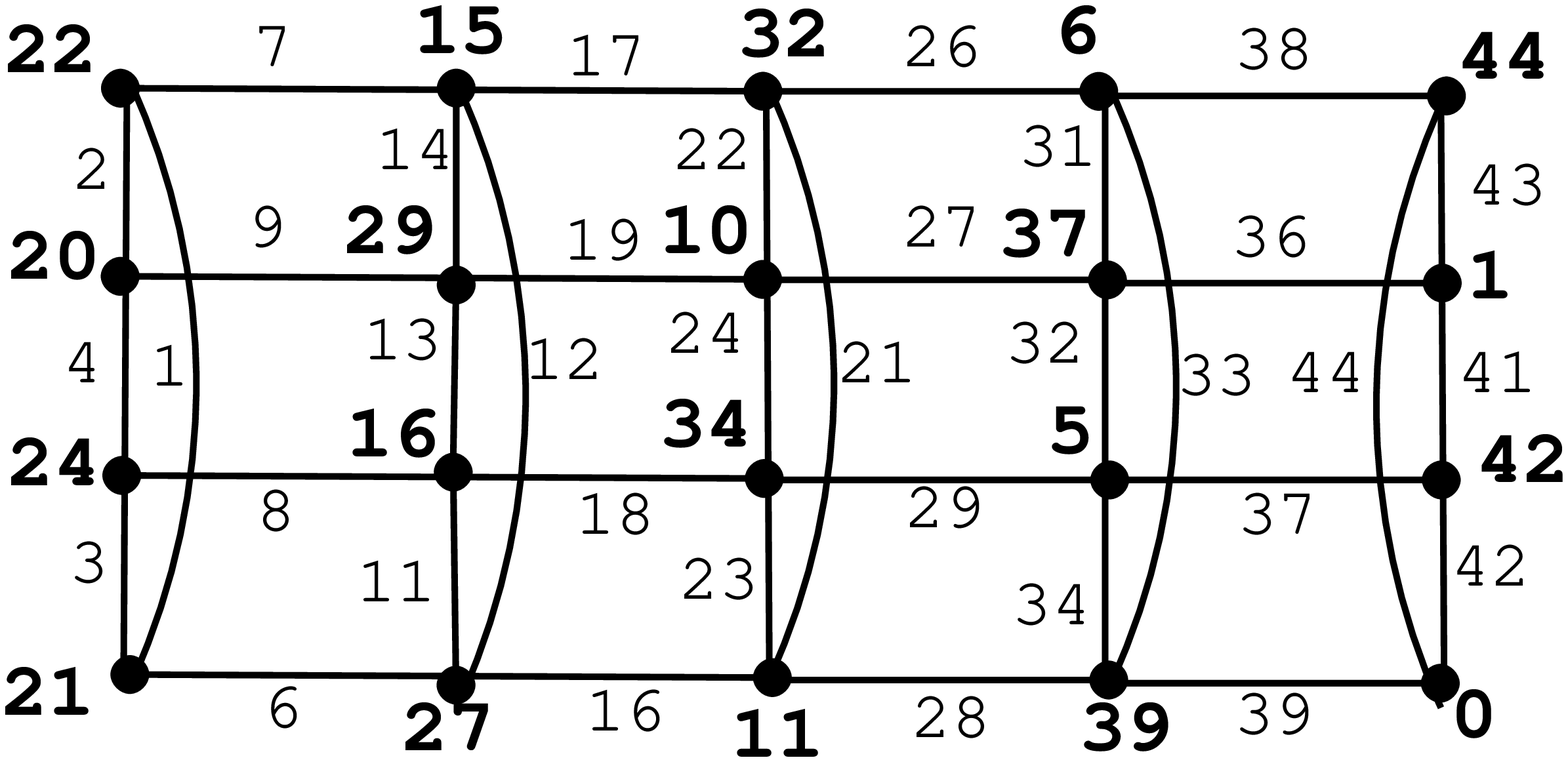}
\caption{}
\label{C4}
\end{center}
\end{figure}
\end{ex}

\begin{thm}\label{d=2(2m-1)}
For any integer $k\geq1$ and $m\geq 2$, $C_{4k}\times P_m$ admits a
$2(2m-1)$-divisible graceful $\alpha$-labeling.
\end{thm}
\noindent
Proof. We will prove the result by induction on $m$.
If $m=2$ the thesis follows from Theorem \ref{prism6}.
Let now $m\geq2$. Suppose that there exists a $2(2m-1)$-divisible graceful $\alpha$-labeling
$f$ of $C_{4k}\times P_m$ with vertices of $C^m$ so labeled:
\begin{align*}
f(C^m) = &(0,(4k+2)(2m-1)-1,1,(4k+2)(2m-1)-2,2,\ldots,\\
&(4k+2)(2m-1)-k,k,(4k+2)(2m-1)-(k+3),
k+1,\dots,\\ &2k-1,(4k+2)(2m-1)-(2k+2)).
\end{align*}
We want to show the existence of a $2(2m+1)$-divisible graceful
$\alpha$-labeling $g$ of $C_{4k}\times P_{m+1}$ satisfying the
same property, namely such that
\begin{eqnarray}
\nonumber g(C^{m+1})&= &(0,(4k+2)(2m+1)-1,1,(4k+2)(2m+1)-2,2,\ldots,\\
\nonumber &\ &(4k+2)(2m+1)-k,k,
(4k+2)(2m+1)-(k+3),k+1,\ldots,\\
 &\ &2k-1,(4k+2)(2m+1)-(2k+2)).\label{2gC(m+1)}
\end{eqnarray}
First of all set
$$g((i,j))=f((i,j))+(4k+2)\quad \forall i=1,\ldots,m,\ \forall j=1,\ldots,4k.$$
This implies that there exists $j\in[1,4n]$ such that $g((m,j))=(4k+2)2m-1$.
We set $g(C^{m+1})$ as in (\ref{2gC(m+1)}) where $g((m+1,j))=0$.
Arguing exactly as in the previous proof one can prove that $g$ is a
$2(2m+1)$-divisible graceful
$\alpha$-labeling $g$ of $C_{4k}\times P_{m+1}$.
\hfill$\Box$

\begin{figure}[H]
\begin{center}
\includegraphics[width=0.4\textwidth]{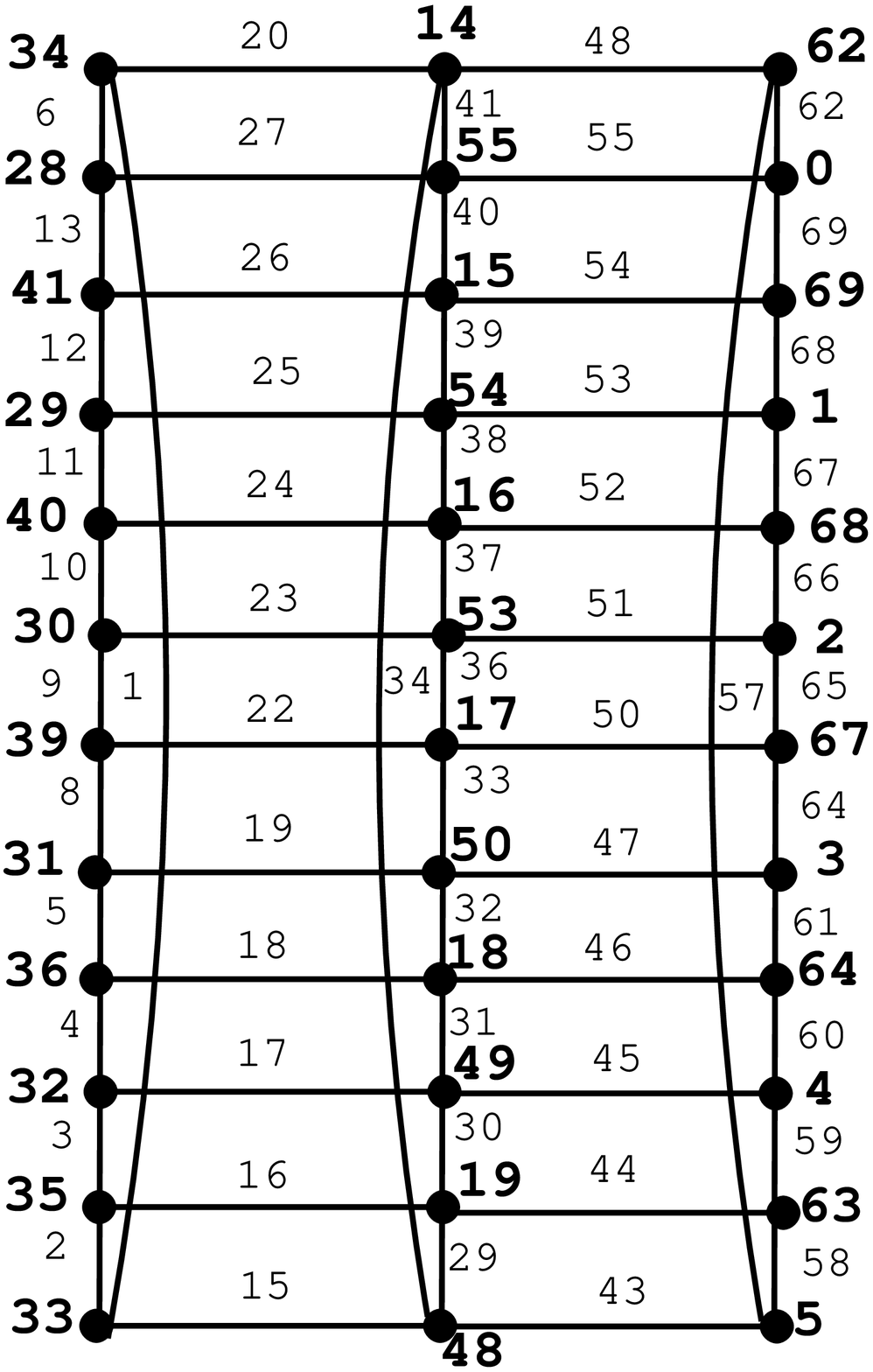}\quad\quad\quad
\includegraphics[width=0.4\textwidth]{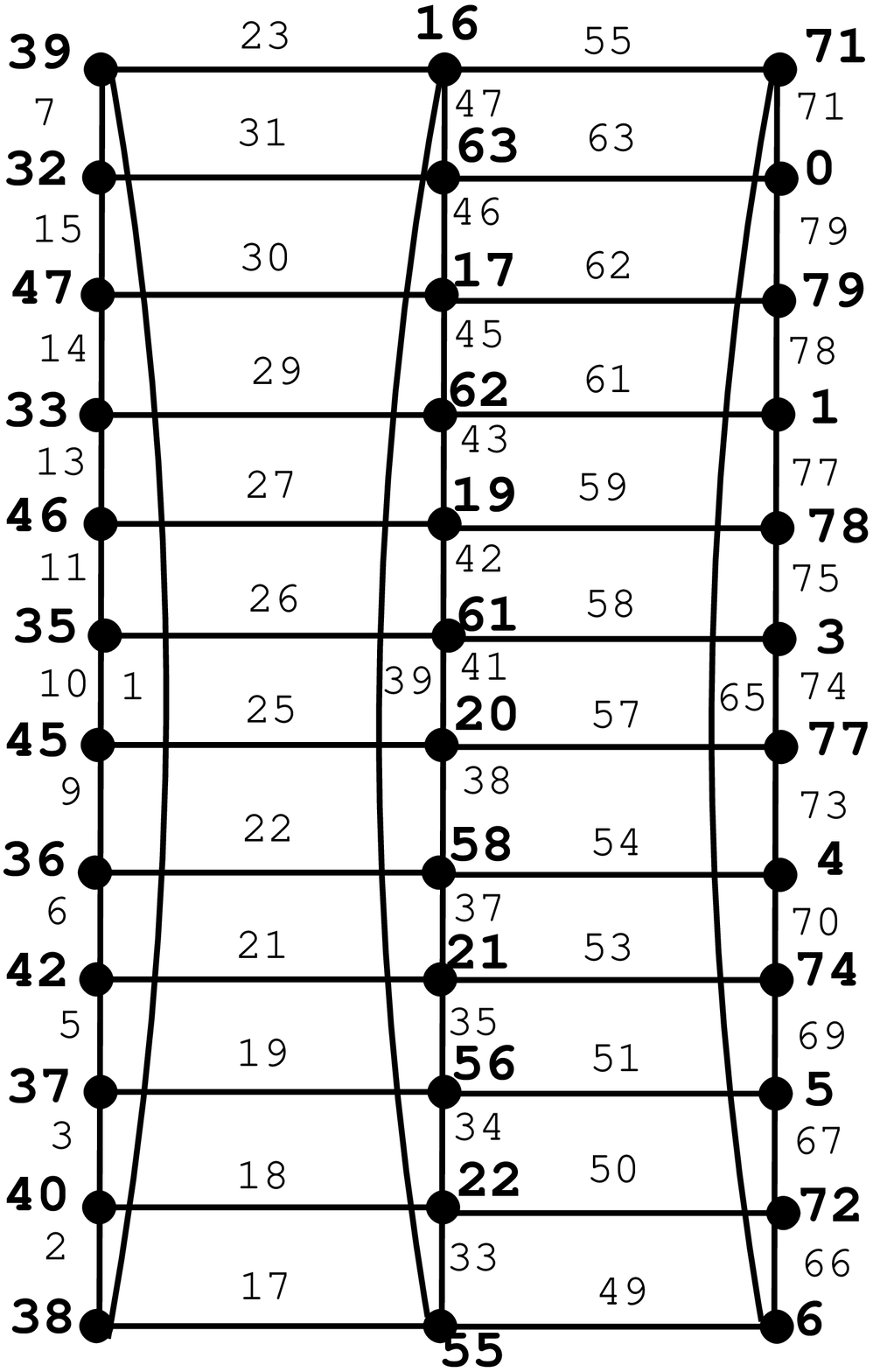}
\caption{A $10$-divisible graceful $\alpha$-labeling of $C_{12}\times P_3$
and a $20$-divisible graceful $\alpha$-labeling of $C_{12}\times P_3$, respectively.}
\label{C12P3}
\end{center}
\end{figure}

\begin{thm}\label{d=4(2m-1)}
For any integer $k\geq1$ and $m\geq 2$, $C_{4k}\times P_m$ admits a
$4(2m-1)$-divisible graceful $\alpha$-labeling.
\end{thm}
\noindent
Proof. We will prove the result by induction on $m$.
If $m=2$ the thesis follows from Theorem \ref{prism12}.
Let now $m\geq2$.
We have to distinguish two cases: $k$ even and $k$ odd.\\
Let $k$ be even.
 Suppose that there exists a $4(2m-1)$-divisible graceful $\alpha$-labeling
$f$ of $C_{4k}\times P_m$ with vertices of $C^m$ so labeled:
\begin{align*}
f(C^m) = &(0,(4k+4)(2m-1)-1,1,(4k+4)(2m-1)-2,2,\ldots,\\
&(4k+4)(2m-1)-\frac{k}{2},\frac{k}{2},(4k+4)(2m-1)-\left(\frac{k}{2}+2\right),
\frac{k}{2}+1,\dots,\\
& (4k+4)(2m-1)-(k+1),k,(4k+4)(2m-1)-(k+4),k+1,\ldots,\\
& (4k+4)(2m-1)-\left(\frac{3}{2}k+2\right),\frac{3}{2}k-1,(4k+4)(2m-1)-\left(\frac{3}{2}k+3\right),\\
& \frac{3}{2}k+1,(4k+4)(2m-1)-\left(\frac{3}{2}k+4\right),\ldots,\\
& (4k+4)(2m-1)-\left(2k+2\right),2k,(4k+4)(2m-1)-(2k+3)).
\end{align*}
We want to show the existence of a $4(2m+1)$-divisible graceful
$\alpha$-labeling $g$ of $C_{4k}\times P_{m+1}$ satisfying the
same property, namely such that
\begin{eqnarray}
\nonumber g(C^{m+1})\hspace{-0.2cm} &= & \hspace{-0.2cm} (0,(4k+4)(2m+1)-1,1,(4k+4)(2m+1)-2,2,\ldots,\\
\nonumber &\ & \hspace{-0.2cm} (4k+4)(2m+1)-\frac{k}{2},\frac{k}{2},(4k+4)(2m+1)-\left(\frac{k}{2}+2\right),
\frac{k}{2}+1,\dots,\\
\nonumber &\ & \hspace{-0.2cm} (4k+4)(2m+1)-(k+1),k,(4k+4)(2m+1)-(k+4),k+1,\ldots,\\
\nonumber &\ & \hspace{-0.2cm} (4k+4)(2m+1)-\left(\frac{3}{2}k+2\right),\frac{3}{2}k-1,(4k+4)(2m+1)-\left(\frac{3}{2}k+3\right),\\
\nonumber &\ & \hspace{-0.2cm} \frac{3}{2}k+1,(4k+4)(2m+1)-\left(\frac{3}{2}k+4\right),\ldots,\\
&\ & \hspace{-0.2cm} (4k+4)(2m+1)-\left(2k+2\right),2k,(4k+4)(2m+1)-(2k+3)).
\label{4gC(m+1)}
\end{eqnarray}
First of all set
$$g((i,j))=f((i,j))+(4k+4)\quad \forall i=1,\ldots,m,\ \forall j=1,\ldots,4k.$$
This implies that there exists $j\in[1,4n]$ such that $g((m,j))=(4k+4)2m-1$.
We set $g(C^{m+1})$ as in (\ref{4gC(m+1)}) where $g((m+1,j))=0$.\\
Let now $k$ be odd. Suppose that there exists a $4(2m-1)$-divisible graceful $\alpha$-labeling
$f$ of $C_{4k}\times P_m$ with vertices of $C^m$ so labeled:
\begin{align*}
f(C^m) = &(0,(4k+4)(2m-1)-1,1,(4k+4)(2m-1)-2,2,\ldots,\\
&(4k+4)(2m-1)-\frac{k-1}{2},\frac{k-1}{2},(4k+4)(2m-1)-\frac{k+1}{2},
\frac{k+3}{2},\\
& \ldots,(4k+4)(2m-1)-k,k+1,(4k+4)(2m-1)-(k+3),k+2,\ldots,\\
& (4k+4)(2m-1)-\frac{3k+3}{2},\frac{3k+1}{2},(4k+4)(2m-1)-\frac{3k+7}{2},\\
& \frac{3k+3}{2},(4k+4)(2m-1)-\frac{3k+9}{2},\ldots,\\
& (4k+4)(2m-1)-\left(2k+2\right),2k,(4k+4)(2m-1)-(2k+3)).
\end{align*}
We want to show the existence of a $4(2m+1)$-divisible graceful
$\alpha$-labeling $g$ of $C_{4k}\times P_{m+1}$ satisfying the
same property, namely such that
\begin{eqnarray}
\nonumber g(C^{m+1})\hspace{-0.3cm} &= &
\hspace{-0.3cm}(0,(4k+4)(2m+1)-1,1,(4k+4)(2m+1)-2,2,\ldots,\\
\nonumber & & \hspace{-0.3cm} (4k+4)(2m+1)-\frac{k-1}{2},\frac{k-1}{2},(4k+4)(2m+1)-\frac{k+1}{2},
\frac{k+3}{2},\\
\nonumber & & \hspace{-0.3cm} \ldots,(4k+4)(2m+1)-k,k+1,(4k+4)(2m+1)-(k+3),k+2,\ldots,\\
\nonumber & & \hspace{-0.3cm} (4k+4)(2m+1)-\frac{3k+3}{2},\frac{3k+1}{2},(4k+4)(2m+1)-\frac{3k+7}{2},\\
\nonumber & & \hspace{-0.3cm} \frac{3k+3}{2},(4k+4)(2m+1)-\frac{3k+9}{2},\ldots,\\
& & \hspace{-0.3cm} (4k+4)(2m+1)-\left(2k+2\right),2k,(4k+4)(2m+1)-(2k+3)).
\label{4gC(m+1)odd}
\end{eqnarray}
First of all set
$$g((i,j))=f((i,j))+(4k+4)\quad \forall i=1,\ldots,m,\ \forall j=1,\ldots,4k.$$
This implies that there exists $j\in[1,4n]$ such that $g((m,j))=(4k+4)2m-1$.
We set $g(C^{m+1})$ as in (\ref{4gC(m+1)odd}) where $g((m+1,j))=0$.\\
Arguing exactly in the proof of Theorem \ref{d=2m-1} one can prove that, in both cases, $g$ is a
$4(2m+1)$-divisible graceful
$\alpha$-labeling $g$ of $C_{4k}\times P_{m+1}$.
\hfill$\Box$\\
\begin{ex}
In Figure \ref{C12P3} we have the $10$-divisible graceful $\alpha$-labeling of $C_{12}\times P_3$ obtained  starting from
the $6$-divisible graceful $\alpha$-labeling of $T_{24}=C_{12}\times P_2$ shown in Figure \ref{T24} and
following the construction
explained in the proof of Theorem $\ref{d=2(2m-1)}$ and
the $20$-divisible graceful $\alpha$-labeling of $C_{12}\times P_3$ obtained  starting from
the $12$-divisible graceful $\alpha$-labeling of $T_{24}$ shown in Figure \ref{T24} and
following the construction
illustrated in the proof of Theorem $\ref{d=4(2m-1)}$.
\end{ex}
\noindent
By virtue of Theorems \ref{decomp}, \ref{d=2m-1}, \ref{d=2(2m-1)} and \ref{d=4(2m-1)}, we have
\begin{prop}
There exists a cyclic  $C_{4k}\times P_m$-decomposition of\\ $K_{(4k+1)\times 2(2m-1)n}$,
 of $K_{(2k+1)\times 4(2m-1)n}$
and of $K_{(k+1)\times 8(2m-1)n}$, for any integers $k,n\geq1$, $m>2$.
\end{prop}

\end{document}